\pdfoutput=1
\documentclass[a4paper,USenglish]{lipics}
\usepackage{microtype}
\usepackage{aslenum}
\usepackage{stmaryrd}

\let\phi\varphi
\newcommand{\syn}[1]{{\color{red}{#1}}}
\newcommand{\sem}[1]{{\color{blue}{#1}}}
\newcommand{\MMM}{\ensuremath{\mathfrak{M}}}
\newcommand{\NNN}{\ensuremath{\mathfrak{N}}}

\title{Syntax versus Semantics}

\titlerunning{Syntax versus Semantics} 

\author[1]{Reinhard Kahle}
\author[2]{Wilfried Keller}

\affil[1]{CMA and DM, FCT, Universidade Nova de Lisboa \\
P-2829-526 Caparica, Portugal \\
\href{mailto:kahle@mat.uc.pt}{kahle@mat.uc.pt}}
\affil[2]{Universität des Saarlandes\\
66123 Saarbrücken, Germany \\
\href{mailto:wilfried.keller@uni-saarland.de}{wilfried.keller@uni-saarland.de}}

\authorrunning{Reinhard Kahle and Wilfried Keller}

\Copyright{Reinhard Kahle and Wilfried Keller}


\serieslogo{logo_ttl}
\volumeinfo
  {M. Antonia {Huertas}, Jo\~ao {Marcos}, Mar\'ia {Manzano}, Sophie {Pinchinat}, \\
  Fran\c{c}ois {Schwarzentruber}}
  {5}
  {4th International Conference on Tools for Teaching Logic}
  {1}
  {1}
  {75}
\EventShortName{TTL2015}

\begin{document}

\maketitle
{\let\thefootnote\relax\footnotetext{\textbf{Acknowledgements}.\enspace The first author was partially supported
by the Fundação para a Ciência e a Tecnologia (Portuguese Foundation for Science
and Technology) through the projects \emph{The Notion of Mathematical Proof}
(PTDC/MHC-FIL/5363/2012) and UID/MAT/00297/2013 (Centro de Matemática
e Aplicações). We thank the anonymous referees for helpful comments.}}

\begin{abstract}
We report on the idea to use colours to distinguish syntax and
semantics as an educational tool in logic classes. This distinction
gives also reason to reflect on some philosophical issues
concerning semantics.
\end{abstract}

\section{Introduction}

When the first author attended as a first year student a logic course
at a philosophy department, he was wondering why the professor was
repeating, after introducing the relation $\vDash$, ``the same'' again
but using the symbol `$\vdash$'. It was, most likely, the negligence of
the student, and not the presentation of the professor, which caused
this fundamental misunderstanding. In retrospect, there 
is one
question which might make one wonder: how was it possible to miss one
of the most fundamental distinctions in modern logic?
On a more mature stage, as third year student, the problem
repeated itself on a different level: speaking about Gödel's first
incompleteness theorem with a teaching assistant, he refused to
continue the discussion when the question was raised what its meaning
is from the perspective of ``proper Mathematics''. Again, the distinction
of syntax and semantics wasn't clearly seen and the discussion stalled
when the participants did not realize that reference of a term like
``proper Mathematics'' need to be better specified while entering into a
discussion of Gödel's theorem.

Based on the experiences above, when starting to teach logic courses
by himself, the first author took it as a particular challenge to
present the problems concerning the distinction of syntax and semantics in a
persuasive way. In this paper, we present the chosen solution:
the use of colours to distinguish the syntactical and semantical role
of logical text. We also report on some of the educational insights
one might win with this approach together with certain philosophical
questions which surface (again) in this context.

We assume that the reader is familiar with first-order logic and has,
at least, an idea of Gödel's incompleteness theorems. Without formal
introduction, we use the standard notations of first-order logic and
Peano Arithmetic, using typically Greek letters from the end of the
alphabet to denote formulae.

\section{First-order logic}

Let us go \emph{in media res} and recall the definition of the
satisfiability relation $\models$ for first-order logic as given by
Barwise in \cite[p.~21]{Bar77}.

\begin{definition}\label{def1}
\renewcommand{\sem}[1]{#1}
\renewcommand{\syn}[1]{#1}
Let \sem{\MMM} be an L-structure.  We define a relation
$$ \MMM \models \phi[s],$$
(read: the assignment $s$ satisfies the formula $\phi$ in $\MMM$) for
all assignments $s$ and all formulae $\phi$ as follows.
\begin{enumerate}[(i)]
\item $\sem{\MMM \models \syn{(t_1=t_2)}[s]}$ iff
  $\sem{\syn{t_1}^{\MMM}(s) = \syn{t_2}^{\MMM}(s)}$,
\item $\sem{\MMM \models \syn{R(t_1,\dots,t_n)}[s]}$,
  iff $\sem{(\syn{t_1}^{\MMM}(s),\dots,\syn{t_n}^{\MMM}(s)) \in \syn{R}^{\MMM}}$,
\item $\sem{\MMM \models \syn{(\neg \phi)}[s]}$  iff
   not $\sem{\MMM \models \syn{\phi}[s]}$,
\item $\sem{\MMM \models \syn{(\phi \wedge \psi)}[s]}$ iff $\sem{\MMM \models \syn{\phi}[s]}$ and $\sem{\MMM \models \syn{\psi}[s]}$,
\item $\sem{\MMM \models \syn{(\phi \vee \psi)}[s]}$ iff $\sem{\MMM
  \models \syn{\phi}[s]}$ or $\sem{\MMM \models \syn{\psi}[s]}$,
\item $\sem{\MMM \models \syn{(\phi \rightarrow \psi)}[s]}$ iff
either not $\sem{\MMM \models \syn{\phi}[s]}$
  or else $\sem{\MMM \models \syn{\psi}[s]}$,
\item $\sem{\MMM \models \syn{(\exists v.\phi)}[s]}$ iff
  there is an $\sem{a\in M}$ such that $\sem{\MMM \models \syn{\phi}[s(^{\sem{a}}_{\syn{v}})]}$,
\item $\sem{\MMM \models \syn{(\forall  v.\phi)}[s]}$ iff for
  all $\sem{a \in M}$, $\sem{\MMM \models \syn{\phi}[s(^{\sem{a}}_{\syn{v}})]}$.
\end{enumerate}
\end{definition}

Barwise then states: ``There is nothing surprising here. It is just
making sure that each of our symbols means what we want it to mean.''

This is certainly correct, but it might be worth to reflect shortly
what this definition actually does.\footnote{Here a disclaimer is
  in order. In the following 
we don't intend to criticise Barwise's
  presentation; first, the citation is taken from a handbook article
  and not from a textbook (and we may add: from an excellent article in an
  excellent handbook); secondly, mathematically there is no issue
  here, the definition is surely the correct one; we are just looking
  for aspects which are of interest in an educational perspective.}
Looking, for instance, to the clauses for conjunction and disjunction,
it simply ``lifts'' the symbols `$\wedge$' and `$\vee$' to the
\emph{meta-level} using the natural language expressions ``and'' and
``or''. Thus, who understands ``and'' and ``or'' is supposed to
understand how to obtain the truth value of a formula using `$\wedge$'
or `$\vee$'. The \emph{meta-level} comes into play, because the natural
language terms are outside of the $\models$ relation, while the
logical symbols are inside.
Nothing surprising here, as said. We will
see later in \S~\ref{sec-int} that the situation for the negation is actually more
subtle than it looks (and this is a first educational pitfall!).
But the issue we like to start with is unambiguously expressed by Barwise
(in continuation of the citation above):
\begin{quote}
There is one possibly confusing point, in (i), caused by our using $=$
for both the real equality (on the right-hand side) and the symbol for
equality (on the left). Many authors abhor this confusion of use and
mention and use something like $\equiv$ or $\approx$ for the symbol.
\end{quote}

It is this confusion\footnote{And it might be doubted, whether the
  difficulty at hand is best couched in terms of use and mention. It
  seems to us rather clear, however, that it is 
  educationally not best dealt with by introducing subtleties of
  quotation theory.}, which will repeat itself when one comes to the
language of arithmetic, we like to address. As Barwise remarks, many
authors try to resolve the problem by the use of different symbols for
syntactical and semantical equality. This would be a way out, but it
comes at the price that the formal language appears artificial, in
particular, if this is extended to arithmetical terms. 
Let us mention the clumsy notation $\textbf{k}_n$ for the numeral
representing the number $n$ in the classical textbook of Shoenfield
\cite{Sho67}. With respect to equality, Sernadas and Sernadas
\cite{SS08}, for instance, use $\cong$ as special sign on the
syntactical side. Cori and Lascar \cite{CL00} use, quite generally, a
bar over the syntactical symbol to refer to its \emph{semantic}
interpretation. 
However, in general (and with good reasons), the syntactic symbols are
chosen to match exactly with their intended meaning; the resulting
``overloading'' is just a problem when the difference of syntax and
semantics is the subject under discussion.

Our solution is to use the same symbols, but
different colours to distinguish the use of symbols on the syntactical
and semantical side. So, let us choose red for syntax and blue for
semantics.\footnote{The particular choice of the colours is arbitrary;
  it is, however, known that one should avoid yellow and green for
  beamer presentations.} After giving a definition of first-order
language with all syntactic expressions in red, and the definition of
a structure with blue for the elements of the structure, the
definition above reads as follows:

\begin{definition}
Let \sem{\MMM} be an L-structure.  We define a relation
$$ \sem{\MMM \models \syn{\phi}[s]},$$
(read: the assignment $\sem{s}$ satisfies the formula $\syn{\phi}$ in $\sem{\MMM}$) for
all assignments $\sem{s}$ and all formulae $\syn{\phi}$ as follows.\footnote{For sure,
    as one referee pointed put, colouring these expressions does not
    remove entirely the
    possibility of misunderstanding: The functional expression
    $\syn{t_1}^{\sem{\MMM}}$, for instance, takes the syntactic argument $\syn{t_1}$ and evaluates it relative to the
    semantic model $\sem{\MMM}$---so the value of the whole expression of course denotes
    a blue item. The interpretation function itself does not appear in this notation, and
    therefore, luckily, we do not have to decide on its colour. However, with interested
    students we had some interesting discussions about its colour.}
\begin{enumerate}[(i)]
\item $\sem{\MMM \models \syn{(t_1=t_2)}[s]}$ iff
  $\sem{\syn{t_1}^{\MMM}(s) = \syn{t_2}^{\MMM}(s)}$,
\item $\sem{\MMM \models \syn{R(t_1,\dots,t_n)}[s]}$,
  iff $\sem{(\syn{t_1}^{\MMM}(s),\dots,\syn{t_n}^{\MMM}(s)) \in \syn{R}^{\MMM}}$,
\item $\sem{\MMM \models \syn{(\neg \phi)}[s]}$  iff
   not $\sem{\MMM \models \syn{\phi}[s]}$,
\item $\sem{\MMM \models \syn{(\phi \wedge \psi)}[s]}$ iff $\sem{\MMM \models \syn{\phi}[s]}$ and $\sem{\MMM \models \syn{\psi}[s]}$,
\item $\sem{\MMM \models \syn{(\phi \vee \psi)}[s]}$ iff $\sem{\MMM
  \models \syn{\phi}[s]}$ or $\sem{\MMM \models \syn{\psi}[s]}$,
\item $\sem{\MMM \models \syn{(\phi \rightarrow \psi)}[s]}$ iff
either not $\sem{\MMM \models \syn{\phi}[s]}$
  or else $\sem{\MMM \models \syn{\psi}[s]}$,
\item $\sem{\MMM \models \syn{(\exists v.\phi)}[s]}$ iff
  there is an $\sem{a\in M}$ such that $\sem{\MMM \models \syn{\phi}[s(^{\sem{a}}_{\syn{v}})]}$,
\item $\sem{\MMM \models \syn{(\forall  v.\phi)}[s]}$ iff for
  all $\sem{a \in M}$, $\sem{\MMM \models \syn{\phi}[s(^{\sem{a}}_{\syn{v}})]}$.
\end{enumerate}
\end{definition}

Now, the difference of $\syn{=}$ and $\sem{=}$ is obvious; the
``possible confusion'' will even catch the student's eye, and this is
of particular educational value. In addition, the quantifier clauses
allow to distinguish better between the syntactical object variable
$\syn{v}$ and the semantical \emph{element} $\sem{a}$---even if
$\sem{a}$ is here, of course, a meta-variable for ``real elements'' of
$\sem{M}$.

Introducing now, as counterpart to the semantic consequence relation,
a derivation relation $\syn{\vdash}$ which
exclusively, deals with red objects, syntactic and semantic reasoning
is already distinguished by colours. It is our teaching experience,
that these colours are of great help for the students, providing some
kind of orientation in rather technical proofs like these for the
compactness and the completeness theorem. As prime examples for the
use of colours we may also mention the Skolem paradox, which, like the
completeness theorem, can be regarded as a corollary of the
compactness theorem.


And so the question of soundness and completeness
is the next one to be addressed.  Even in a compact course on Gödel's
incompleteness theorems these results of Gödel's dissertation should
be discussed.\footnote{In introductory logic courses in philosophy
  seminars these results will usually not be proven, but, of course,
  deserve at least to be mentioned and discussed.} Today, they are
usually not proven along the lines of Gödel's original work, but
rather with the well-known technique introduced by Henkin. This
technique is standard for many completeness proofs for quite a couple
of logics and, therefore, constitute an integral tool in the toolbox
of the mathematical logician: it proceeds by constructing a blue world
(validating the negation of a formula, that cannot be proven in the
syntactic calculus) out of the red material itself.
(The situation is somewhat dual to the more demanding case of (set-theoretic)
forcing, where one, so to speak, bestows properties of a language on certain
blue sets, thereby emulating the red world with originally blue objects.)

\section{Arithmetic} 

The use of the colours exhibits its real potential only, when we study
Arithmetic and want to motivate, for instance, Gödel's incompleteness
theorems for Peano Arithmetic. 
In its original form, Gödel's proof actually entirely refrains from
the semantic realm and in some sense requires only
``red'' reasoning.
It, indeed, is a theorem about the red
  world, but one that operates on a meta-level. It speaks about the
  axiomatic system, but it does not reason about its meaning, quite
  like Hilbert's attempted consistency proofs do. And likewise it uses
  only finitistically acceptable means.\label{Hilbert}
In \cite[fn.~1, p.~337]{Bea02} the editors 
explicitly
attribute Gödel's avoiding of the concept of truth in his 1931 paper
to the ``respect for the `prejudice' of the Hilbert school.'' Further
reflections 
can be found in \cite{Fef84}.

But, of course, the first incompleteness theorem
can---and probably should---be motivated by the question whether Peano
Arithmetic is complete with respect to the \emph{standard model}.
Peano Arithmetic may be defined over the set
$\{\syn{0},\syn{S},\syn{+},\syn{\cdot}\}$ of non-logical symbols.
The standard model, as semantic object, can be given in blue as
follows:
$\sem{\NNN = \langle \mathbb{N},0,S,+,\cdot\rangle}.$
It goes without saying that the red symbols should be interpreted by
their respective blue counterparts.
Gödel's first incompleteness theorem tells us, that the red world, as
long as a consistent recursive axiomatic system is chosen, will miss
some formula \emph{and} its negation\footnote{Of course on pain of
  inconsistency it should \emph{not} prove any formula together with
  its negation.}, 
it 
cannot decide 
all given formulae by proof. However, the blue world
so to speak decides formulae: Either a formula or its negation is
satisfied (in a model by a valuation); this is a direct
  consequence of clause  (iii) in Def.\ \ref{def1} (see also 
\S\ref{sec-int} below).
There is no gap in the blue world---but there will always be one in the
(sensibly chosen, i.e., recursive and consistent) red one.
This then not only entails that no red world (as above) fixes the blue
world of Arithmetic as its only model---
this could already be concluded from
Löwenheim--Skolem---but that not all red statements that are true
about the blue world can be proven in the red realm.

So, \emph{how} are the blue objects in $\sem{\NNN = \langle
  \mathbb{N},0,S,+,\cdot\rangle}$ actually given? $\sem{\mathbb{N}}$
  as the set of natural numbers $\sem{\{0,1,2,\dots\}}$ can
be taken for granted, and with it, of course, also the object
$\sem{0}$. But what \emph{is} $\sem{+}$? It should be, of course, the
set of all triples $(\sem{a},\sem{b},\sem{c})$ such that $\sem{a} +
\sem{b} = \sem{c}$. But this doesn't help us, as we just moved the
question to another level, writing $+$ in black (this plus
  sign can, of course, not be blue if we don't want to run into an
  immediate circle). To avoid it, one could get back to use of dots
and say that $\sem{+}$ is the set:
\begin{align*}
\{&\sem{(0,0,0)}, \sem{(0,1,1)}, \sem{(0,2,2)}, \sem{(0,3,3)}, \dots
\\
& \sem{(1,0,1)}, \sem{(1,1,2)}, \sem{(1,2,3)}, \sem{(1,3,4)}, \dots \\
& \sem{(2,0,2)}, \sem{(2,1,3)}, \sem{(2,2,4)}, \sem{(2,3,5)}, \dots \\
& \qquad \vdots \hspace*{150pt} \}.
\end{align*}

In fact, we see only two ways out of this situation: First, we appeal
to the common platonism in mathematics and stipulate that $\sem{+}$
exists in the structure $\sem{\NNN}$ in the way we expect it and with
the properties we want it to have.  Or, secondly, we treat $\sem{+}$ as a
proper set-theoretic object and construct it within our favorite
axiomatic set theory, let's say $\textsf{ZFC}$.

While the second one is a solution concerning the level of Arithmetic, it simply
shifts the problem to set theory. Within $\textsf{ZFC}$, the defined $\sem{+}$
would be a red object (relative to $\textsf{ZFC}$), and we would have to ask now,
how the models of $\textsf{ZFC}$ look like---which is, not surprising, a much
harder question than to look for structures for Arithmetic, not only
in a technical sense, but also in a philosophical one; in particular,
it leads us simply back to the situation we are discussing for
$\sem{+}$ just for other, now set-theoretical, objects.

For the first option, the sheer existence of $\sem{+}$ presupposes a
certain form of platonism. This platonism allows to dispose of any
further technical questions, but one enters a dangerous
philosophical terrain.

We may leave aside here the preferred way out (this is a philosophical
question---but one which is naturally asked in this context and can
nicely be brought out and talked about with the help of our colours);
what we consider as an important lesson for the teaching of
logic is, that any discussion of the syntax/semantics distinction in
the context of Arithmetic---and, \emph{a fortiori}, in relation to
Gödel's Incompleteness Theorems---has to presuppose some ``blue
objects''. And this presupposition is far from being trivial. As
a minimal conclusion, we like to state that an unreflected ``identification'' of,
for instance, the ``syntactical addition'' (symbol) $\syn{+}$ with the
``semantical addition'' (function) $\sem{+}$ prevents a student from
any understanding of the problem which is at issue here.

\section{Gödel's Incompleteness Theorems}

Gödel's Incompleteness Theorems are, without doubt, the most important
results in mathematical logic. They are, however, clouded by continuing
philosophical misunderstandings,\footnote{For a detailed discussion,
  in particular in view of Gentzen's later results, we refer to
  \cite{Kah15}.} and a pure formal presentation of the
results can easily support such prejudices.\footnote{An extreme case
  was reported from a philosophical seminar: the participants,
  after studying carefully---and probably painfully---a formal
  presentation of Gödel's proof, came to the conclusion that, yes,
  after all the results seem to be
  correct \emph{if one has a language with 7 symbols; but, of course,
    one doesn't know whether it still holds when one would have 8
    symbols}.}

We like to address only one problem coming from the fact that many
presentations make use of the \emph{Chinese remainder theorem}.
Apparently, here proper mathematics is used to prove something about
formal systems for mathematics and one may wonder whether there is
some kind of vicious circle in the back.\footnote{The first author
  remembers an elaboration of a Math student of a course presentation
  of Gödel's theorems which focussed nearly entirely on the use of the
  Chines remainder theorem; apparently it was the only part of the
  proof the student was properly understanding based on his Math
  classes.}  This is, of course, not the case; in fact, the use of
this theorem is only relevant for a certain technical step, in
particular in the case of Peano Arithmetic, and it is irrelevant for
other axiomatization or codings (see \cite[\S3.2.6]{Smo77}). Hao Wang reported
this from interviews with Gödel \cite[p.~653]{Wan81}: ``He enjoyed
much the lectures by Furtw[ä]ngler on number theory and developed an
interest in this subject which was, for example, relevant to his
application of the Chinese remainder theorem in expressing primitive
recursive functions in terms of addition and multiplication.'' Thus,
what is at issue here is the representation of primitive recursive
functions in the formal arithmetical theory. Albeit, it is, of course,
possible to develop primitive recursion \emph{in} such theories, a
proper educational approach \emph{can} ``outsource'' primitive recursion.

We did that by setting up a new realm for primitive recursion
(formally a functional algebra, consisting of function symbols and
equalities)---and using a new colour for it, as, for the moment, it
can be located outside of the red and blue world. After providing a
sufficient stock of primitive recursive functions, one simply proves
the representation theorem showing that every such function can be
represented in the formal theory under consideration (as Peano
Arithmetic, for instance). In the particular case of Peano Arithmetic,
this proof may use the Chinese remainder theorem---but it is obvious
that it is not relevant if one would provide another proof, in
particular for other theories. In addition, it should become clear
that Gödel's theorems can be carried out in essentially the same way for
all theories allowing for the appropriate representation theorem.

With respect to the Chinese Remainder Theorem and the like the
diagnosis is as following: The Chinese Remainder Theorem is---as far
as Gödel's first incompleteness theorem is concerned---at the level of
meta theory, but not necessarily in the blue world. It
deals with the
recursion-theoretic realm, 
concerning codings, and can be
outsourced and imported if needed.

The situation changes however, if one wants to proceed to the second
Incompleteness Theorem: There one repeats the proof of the first
completely inside the red theory. For that case it is of vital
importance that one uses codings and theorems that can be reproduced
in the object theory itself---this is the sense in which it is
``occasionally not just important [in logic] what you prove, but how
you prove it'' \cite[p.~190]{rautenberg2006concise}.

\section{``Semantics comes first''?}

The dictum ``Semantics comes first'' is attributed to
Tarski.\footnote{The well-known computer software package
  \emph{Tarski's World} \cite{BPBE08} illustrates this slogan quite
  appropriately when the student can literally manipulate the
  semantical objects.} And it has some rationale: in first-order logic
as well as in Arithmetic, we presuppose the \emph{meanings} of the
logical and arithmetical expressions before we start to manipulate
them. It is, in general, seen as a task of the axiomatization to
provide calculi which capture such meaning---and the completeness
theorems should show us that the axiomatization was successful.

As much as completeness theorems are concerned, they give us---\emph{a
  posteriori}---the possibility to identify the red with the blue
objects. But, to make sense out of a completeness theorem---and to
prove it---the difference of the syntactical and semantical expression
has to be seen clearly. In this perspective, we could even propose
that one could forget of the red and blue colour for, let say,
equality, \emph{after} the proof of the completeness theorem for
first-order logic.


But there is still an issues to discuss which might call Tarski's
dictum in question.

First, the dictum makes sense only, if the semantics is
\emph{understood before} the syntax. It is doubtful that
semantics could serve its purpose if it is, by itself, not fully
understood. For first-order logic and Arithmetic it should be the case
that their semantics is clear (although we mentioned the problem to
precisely articulate it in the case of $\sem{+}$, for instance). And this holds
probably also for most of the mathematical theories we are familiar
with.

The situation changed drastically, when Computer Science entered the
stage. Semantics of programming languages is a rather
challenging topic in CS. And here, the syntax is essentially always
prior to the semantics. Instead of finding an axiomatization for a given
semantics, now one looks for a semantics for a given programming
language. But it is not only that such a semantics is no longer
``first'', it might be the case that proposed semantics
are, in a technical sense, more complicated than the programming
language in itself (as example, we may refer to
  \cite{Aea11}). In some cases, one might wonder in which sense
such semantics provide any \emph{meaning}.\footnote{A colleague
  expressed this in a sarcastic statement: ``They put these
  brackets $\llbracket \cdot \rrbracket$ around something and believe they
  have achieved something.''}

From the educational perspective, this only means that logic cannot
any longer be taught with Tarski's dictum as a rock solid starting
point. The differentiation of syntax and semantics becomes even more
significant. What it means to put ``syntax first'' we will discuss in
the next section.

\section{Intuitionism}\label{sec-int}

When Brouwer conceived his mathematical philosophy of
\emph{intuitionism} he had definitely nothing formal in mind. However,
based on the ``successful'' axiomatization of intuitionistic logic by
Heyting, today, intuitionism is often presented in an axiomatic
context; for propositional logic, for instance, just like a
  calculus for classical logic, but without \emph{tertium-non-datur}.
In this way, it can be recast totally in our red world. An
``additional'' semantics, may it be informal or set-theoretically, would
be alien to intuitionism.\footnote{The first author remembers a remark
  of an author to the effect that intuitionism was not succeeding due to the lack of
  a proper semantics, until this deficit was cured by the introduction
  of Kripke semantics. In our view, this sounds like claiming that a
  protestant church is not a church as long as it doesn't have a
  pope.} One could go a step further and say---instead of that there
is no blue world---that the red world should be its own blue
world. This is often paraphrased by saying that, from an intuitionistic
perspective, the meaning of a formula is given by its proof (better:
the set of all its proofs).
Today, this idea is reflected in \emph{proof-theoretic semantics}, an
approach which ``attempts to locate the meaning of propositions and
logical connectives not in terms of interpretations \dots  but in the
role that the proposition or logical connective plays within the
system of
inference.''\footnote{\texttt{http://en.wikipedia.org/wiki/Proof-theoretic\_semantics},
  accessed Febr.\ 20, 2015. For more on proof-theoretic semantics, see
\cite{KSH06,SH14}.}

In view of our discussion above, the formal framework of intuitionism
undermines substantially Tarski's dictum. In particular, when we go
back to Definition~\ref{def1} of a L-structure, one can easily observe
that clause (iii) for negation \emph{builds in} classical logic into
any structure.

It is worth to reflect a little bit more on this clause; it
looks as innocent as any other clause, that for conjunction and
disjunction, for instance,
``just making sure that each of our symbols means what we want it to
mean.'', to repeat the citation from Barwise \cite[p.~21]{Bar77}. But
it is not only the case, that Brouwer would not agree that this is what
we want negation to mean; it is the problem that the right-hand side---not
$\sem{\MMM \models \syn{\phi}[s]}$---can, in general, not practically
be verified. It is the full purpose of an axiomatization to give a
\emph{positive} approach to negation, i.e., the symbol for negation
$\syn{\neg}$ has to stay always on the right-hand side of the
derivation sign $\syn{\vdash}$, and we must not make use of a
``meta-negation'' $\syn{\nvdash}$. That this is possible is far from
being trivial---the complexity is visible in the proof of Gödel's
completeness theorem---but it makes clear the fundamental difference
between syntax (in form of an axiomatization, providing
an---executable, albeit not decidable---derivation relation) and
semantics, which is, in general, inherently non-constructive.

As much as it comes to Brouwer's criticism of the usual semantic
understanding of negation, it is worth to cite Bernays \cite[p.~4 (our
translation)]{Ber79}, who replied to
it as follows:
\begin{quote}
As one knows, the use of the ``tertium-non-datur'' in relation to
infinite sets, in particular in Arithmetic, was disputed by
\textsc{L.\ E.\ J.\ Brouwer}, namely in the form of an opposition to
the traditional logical principle of the excluded middle. Against this
opposition it is to say that it is just based on a reinterpretation of
the negation. \textsc{Brouwer} avoids the usual negation non-A, and
takes instead ``A is absurd''. It is then obvious that the general
alternative ``Every sentence A is true or absurd'' is not justified.
\end{quote}

Following Bernays, the controversy about classical and intuitionistic
negation could even be boiled down to a question of words. What is of
importance for us here, is that the semantic determination of
(classical) negation closes the road to understand Brouwer's
(intuitionistic) negation. To say it differently: the ``standard
approach'' of teaching logic based on Tarski's dictum (used for
Tarskian semantics as given in Def.\ \ref{def1}) blocks a proper
transition to intuitionistic logic, at least, if the students are
``semantically biased''. 
As a matter of fact, intuitionistic logic can---and maybe: has to---be
developed purely in the red world.\footnote{As said at the beginning,
  this concerns intuitionistic logic based on Heyting's
  axiomatization; Brouwer's original intuitionism, in particular if
  considered in the mathematical area of Analysis, should not come in
  a formal livery at
  all.}
It is not our aim to motivate an approach that exclusively
focusses on intuitionism or the like.\footnote{It is worth mentioning
  that \emph{dialogical logic} has a rather interesting status
  concerning the syntax and semantics distinction, which doesn't seem to be
  properly explored yet, see Kahle \cite{Kah08,Kah15l}.}
But it is important that an
introduction to classical logic should not come with a semantical bias
which blocks the road to non-classical logic from the
start.\footnote{For instance,
substructural logics are best approached via proof theory, because---quite
like intuitionism---there lies their motivation; on the other hand their
semantics usually turn out to be rather complicated.}
Here the colours may help.

\section{Conclusion}

It is our experience from courses on mathematical logic---regular
university classes, block courses, and summer school courses---that the
use of colours gives the students, indeed, a ``tool'' at hand to see
better the fundamental distinction of syntax and
semantics.\footnote{We have, of course, not developed an empirical study
  to ``prove'' that the use of colours improves the teaching; it is
  rather an informed impression backed by feedback from our students.
  And although we do not have statistically valid data comparing the
  effects of this method with other approaches, the test results of the
  courses were very positive (but we had few, but rather good
  and understanding students).}
Although the use of colours requires some additional effort in the
preparation of the student material, it is a rather inexpensive
investment to obtain significant educational added value. At a certain
stage---for instance, as mentioned, after the proof of the
completeness theorem---, when the student has already internalized the
difference of syntax and semantics one can even easily go back to the
common identification of the syntactical and semantical expression
using them just in black; it provides even a good test template, as
one can ask the students in all kind of instances whether a
certain expression would have to be coloured red or blue.

In addition, while the colouring can, of course, not solve any
philosophical problem concerning syntax and semantics, it helps to
sharpen the sensibility for the underlying problems. In any case,
the success of a logic course will not only depend on the use of tools,
but on the interaction between teacher, tutor(s) and students.
The use of a new dimension---like colour---underlines the importance
of the distinction between syntax and semantics, it brings various
(philosophical) topics into focus and gives a convenient way of
talking about these issues---a way that often is readily accepted
by the students.

The message of this paper is not: By
all means use colours to distinguish syntax from semantics!
But rather:
Distinguish syntax and semantics clearly and discuss some of the
issues related to that distinction---the use of colour will help you
doing so!


\bibliography{svs}

\begin{thebibliography}{10}

\bibitem{Aea11}
Th. Altenkirch, P.~Morris, F.~Nordvall Forsberg, and A.~Setzer.
\newblock A categorical semantics for inductive-inductive definitions.
\newblock In A.~Corradini and B.~Klin, editors, {\em Algebra and Coalgebra in
  Computer Science}, volume 6859 of {\em Lecture Notes in Computer Science},
  pages 70--84. Springer, 2011.

\bibitem{BPBE08}
D.~Barker-Plummer, J.~Barwise, and J.~Etchemendy.
\newblock {\em Tarski's World}, volume 169 of {\em CSLI Lecture Notes}.
\newblock CSLI Publications, 2008.

\bibitem{Bar77}
J.~Barwise.
\newblock An introduction to first-order logic.
\newblock In J.~Barwise, editor, {\em Handbook of Mathematical Logic}, pages
  5--46. North-Holland, 1977.

\bibitem{Ber79}
P.~Bernays.
\newblock {Bemerkungen zu \textsc{Lorenzen}'s Stellungnahme in der Philosophie
  der Mathematik}.
\newblock In K.~Lorenz, editor, {\em Konstruktionen versus Positionen},
  volume~1, pages 3--16. Berlin, 1979.

\bibitem{CL00}
R.~Cori and D.~Lascar.
\newblock {\em Mathematical Logic, Part I}.
\newblock Oxford University Press, 2000.

\bibitem{Bea02}
B.~Buldt et~al., editor.
\newblock {\em {Kurt G{\"o}del. Wahrheit und Beweisbarkeit}}, volume 2:
  Kompendium zum Werk.
\newblock {\"obv \& hpt}, 2002.

\bibitem{Fef84}
S.~Feferman.
\newblock {Kurt G{\"o}del: Conviction and Caution}.
\newblock {\em Philosophia Naturalis}, 21:546--562, 1984.

\bibitem{Kah08}
R.~Kahle.
\newblock {Konstuktivismus und Semantik}.
\newblock In J.~Mittelstra{\ss}, editor, {\em {Der Konstruktivismus im Ausgang
  der Philosophie von Wilhelm Kamlah und Paul Lorenzen}}, pages 197--212.
  Mentis, 2008.

\bibitem{Kah15l}
R.~Kahle.
\newblock {Dialoge als Semantik}.
\newblock In J.~Mittelstra{\ss} and Chr. von B{\"u}low, editors, {\em
  {Dialogische Logik}}, pages 43--54. Mentis, 2015.

\bibitem{Kah15}
R.~Kahle.
\newblock Gentzen's consistency proof in context.
\newblock In R.~Kahle and M.~Rathjen, editors, {\em Gentzen's Centenary}.
  Springer, 2015, to appear.

\bibitem{KSH06}
R.~Kahle and P.~Schroeder-Heister, editors.
\newblock {\em Proof-theoretic semantics}, volume 148.
\newblock 2006.

\bibitem{rautenberg2006concise}
W.~Rautenberg.
\newblock {\em A Concise Introduction to Mathematical Logic}.
\newblock Springer, 2006.

\bibitem{SH14}
P.~Schroeder-Heister.
\newblock Proof-theoretic semantics.
\newblock In E.\~N.\ Zalta, editor, {\em The Stanford Encyclopedia of
  Philosophy}. Summer 2014 edition, 2014.

\bibitem{SS08}
A.~Sernadas and C.~Sernadas.
\newblock {\em Foundations of Logic and Theory of Computation}, volume~10 of
  {\em Texts in Computer Science}.
\newblock College Publication, 2008.

\bibitem{Sho67}
J.~R. Shoenfield.
\newblock {\em Mathematical Logic}.
\newblock Addison-Wesley, Reading, MA, 1967.
\newblock Reprinted: Association for Symbolic Logic and AK Peters, 2001.

\bibitem{Smo77}
C.~Smorynski.
\newblock The incompleteness theorems.
\newblock In J.~Barwise, editor, {\em Handbook of Mathematical Logic}, pages
  821--865. North-Holland, 1977.

\bibitem{Wan81}
H.~Wang.
\newblock Some facts about {K}urt {G}{\"o}del.
\newblock {\em The Journal of Symbolic Logic}, 46:653--659, 1981.

\end{thebibliography}

\newpage
\thispagestyle{empty}
{\ }

\end{document}